\pdfoutput=1
\RequirePackage{ifpdf}
\ifpdf 
\documentclass[pdftex]{sigma}
\else
\documentclass{sigma}
\fi

\numberwithin{equation}{section}

\newtheorem{Theorem}{Theorem}[section]

\newtheorem{Proposition}[Theorem]{Proposition}
\newtheorem{Conjecture}[Theorem]{Conjecture}
 { \theoremstyle{definition}

\newtheorem{Remark}[Theorem]{Remark} }

\usepackage{stmaryrd}
\usepackage[mathscr]{eucal}
\usepackage{calligra}

\newcommand{\tr}{\operatorname{tr}}

\newcommand{\pr}{p}

\newcommand{\Id}{\operatorname{Id}}

\newcommand{\Hom}{\operatorname{Hom}}

\newcommand{\ad}{\operatorname{ad}}

\newcommand{\dbar}{\bar{\partial}}

\newcommand{\CC}{{\mathbb C}}

\newcommand{\PP}{{\mathbb P}}

\newcommand{\RR}{{\mathbb R}}

\newcommand{\rk}{\operatorname{rk}}

\renewcommand{\)}{\right)}

\newcommand{\vol}{\operatorname{vol}}

\newcommand{\Vol}{\operatorname{Vol}}

\newcommand{\defeq}{\mathrel{\mathop:}=} 

\newcommand{\surj}{\to\kern-1.8ex\to}

\newcommand{\lto}{\longrightarrow}

\newcommand{\lra}[1]{\stackrel{#1}{\longrightarrow}}

\newcommand{\cO}{\mathcal{O}}

\newcommand{\Lie}{\operatorname{Lie}}

\newcommand{\LieG}{\operatorname{Lie} \GGG}

\newcommand{\LieX}{\operatorname{Lie} \widetilde{\GGG}}

\newcommand{\LieH}{\Lie\HHH}

\newcommand{\SL}{\operatorname{SL}}

\newcommand{\U}{\operatorname{U}}

\newcommand{\SU}{\operatorname{SU}}

\DeclareFontFamily{OT1}{rsfs}{}
 \DeclareFontShape{OT1}{rsfs}{n}{it}{<->rsfs10}{}
 \DeclareMathAlphabet{\curly}{OT1}{rsfs}{n}{it}

\newcommand{\AAA}{\curly{A}}
\newcommand{\JJJ}{\curly{J}}

\newcommand{\GGG}{\curly{G}}
\newcommand{\HHH}{\curly{H}}
\newcommand{\III}{\curly{I}}
\newcommand{\PPP}{\curly{P}}

\begin{document}

\allowdisplaybreaks

\renewcommand{\thefootnote}{}

\newcommand{\arXivNumber}{2309.15673}

\renewcommand{\PaperNumber}{032}

\FirstPageHeading

\ShortArticleName{K\"ahler--Yang--Mills Equations and Vortices}

\ArticleName{K\"ahler--Yang--Mills Equations and Vortices\footnote{This paper is a~contribution to the Special Issue on Differential Geometry Inspired by Mathematical Physics in honor of Jean-Pierre Bourguignon for his 75th birthday. The~full collection is available at \href{https://www.emis.de/journals/SIGMA/Bourguignon.html}{https://www.emis.de/journals/SIGMA/Bourguignon.html}}}

\Author{Oscar GARC\'IA-PRADA}

\AuthorNameForHeading{O. Garc\'ia-Prada}

\Address{Instituto de Ciencias Matem\'aticas (CSIC-UAM-UC3M-UCM),\\
Nicol\'as Cabrera 13--15, Cantoblanco, 28049 Madrid, Spain}

\Email{\href{mailto:oscar.garcia-prada@icmat.es}{oscar.garcia-prada@icmat.es}}
\URLaddress{\url{https://www.icmat.es/miembros/garcia-prada/}}

\ArticleDates{Received October 02, 2023, in final form April 04, 2024; Published online April 11, 2024}

\Abstract{The K\"ahler--Yang--Mills equations are coupled equations for a K\"ahler metric on a compact complex manifold and a connection on a complex vector bundle over it. After briefly reviewing the main aspects of the geometry of the K\"ahler--Yang--Mills equations, we consider dimensional reductions of the equations related to vortices --- solutions to certain Yang--Mills--Higgs equations.}

\Keywords{K\"ahler--Yang--Mills equations; vortices; gravitating vortices; dimensional reduction; stability}

\Classification{32Q20; 53C07}

\begin{flushright}
\begin{minipage}{65mm}
\it Dedicated to Jean-Pierre Bourguignon\\
 with gratitude and admiration
 \end{minipage}
\end{flushright}

\renewcommand{\thefootnote}{\arabic{footnote}}
\setcounter{footnote}{0}

\section{Introduction}\label{sec:intro}

I would like to single out two contributions of Jean-Pierre Bourguignon related to this paper. The first one is concerned with the study of properties of the scalar curvature of a Riemannian metric \cite{bourguignon-ezin}, more concretely, with the problem of prescribing the scalar curvature in a conformal class (what is referred as the {\em Nirenberg problem} in the case of the
standard conformal class on the 2-sphere); the second contribution regards the study of Yang--Mills connections \cite{bourguignon-lawson}.
In this paper we consider a system of partial differential equations that somehow combines these two problems in a K\"ahlerian
set up.
These equations, known as the {\em K\"ahler--Yang--Mills equations}, were introduced in \cite{AGG} based on the PhD Thesis of M. Garc\'ia-Fern\'andez (Madrid, 2009) \cite{GF}. The K\"ahler--Yang--Mills equations are coupled equations for a K\"ahler
metric on a compact complex manifold and a connection on a complex vector bundle over it.
They emerge from a natural extension of the theories for constant scalar curvature K\"ahler metrics (Yau--Tian--Donaldson)
and Hermite--Yang--Mills connections (Donaldson--Uhlenbeck--Yau).

Fix a holomorphic vector bundle $E$ over a compact complex K\"ahlerian manifold $M$. The K\"ahler--Yang--Mills equations intertwine the scalar curvature $S_g$ of a K\"ahler metric $g$ on $M$ and the curvature $F_H$ of the Chern connection of a Hermitian
metric $H$ on $E$:
\begin{equation}\label{eq:CKYM0}
{\rm i}\Lambda_g F_H = \lambda \Id,\qquad
S_g - \alpha \Lambda^2 \tr F_H^2 = c.
\end{equation}
Here, $\Lambda_g F_H$ is the contraction of the curvature $F_H$ with the K\"ahler form of $g$. The equations depend on a coupling constant $\alpha\in\RR$, and the constants $\lambda, c\in\RR$ are topological.

Equations \eqref{eq:CKYM0} have a symplectic interpretation. As
observed by Fujiki~\cite{Fujiki} and Donaldson~\cite{D1}, the constant scalar curvature condition for a
K\"ahler metric has a moment map interpretation in terms of a symplectic form
$\omega$ on the smooth compact manifold $M$. The group of symmetries
of the theory for constant scalar curvature
K\"ahler metrics is the group $\HHH$ of {\em Hamiltonian
symplectomorphisms} of $(M,\omega)$. This group acts on the space $\JJJ^i$ of {\em integrable
almost complex structures} on $M$ which are compatible with $\omega$,
and this action is Hamiltonian for a natural symplectic form
$\omega_\JJJ$ on $\JJJ^i$. The moment map interpretation of the Hermite--Yang--Mills
equation was pointed out first by Atiyah and Bott~\cite{AB} for the
case of Riemann surfaces and generalized by Donaldson~\cite{D3} to
higher dimensions. Here, one considers the symplectic action of the
{\em gauge group} $\GGG$ of the Hermitian bundle $(E,h)$ on the space of {\em unitary connections}~$\AAA$
endowed with a natural symplectic form $\omega_\AAA$. Relying on these
two cases, the phase space for the K\"ahler--Yang--Mills theory is provided by the
subspace of the product $\PPP \subset \JJJ^i \times \AAA$
defined by the additional integrability condition for a connection $A\in \AAA$ given by the vanishing of the~$(0,2)$-part of its
curvature. Our
choice of symplectic structure is the restriction to $\PPP$ of the
symplectic form \smash{$\omega_\alpha = \omega_\JJJ + \frac{4}{(n-1)!}\alpha \omega_\AAA$},
for a non-zero coupling constant $\alpha\in \RR$. Here $n$ is the complex dimension of~$M$.

Consider now the {\em extended gauge group} $\widetilde{\GGG}$,
defined as the group of automorphisms of the Hermitian bundle $(E,H)$ covering
Hamiltonian symplectomorphisms of $(M,\omega)$. This is a non-trivial
extension
\begin{equation}\label{eq:Ext-Lie-groupsintro}
1 \to \GGG \to \widetilde{\GGG} \to \HHH \to 1,
\end{equation}
where $\GGG$ is the group of automorphisms of $(E,H)$ covering the identity
on $M$ --- the usual gauge group ---, and $\HHH$, as above, is the group of Hamiltonian
symplectomorphisms of $(M,\omega)$. The group $\widetilde{\GGG}$ acts on
$\PPP$ in a Hamiltonian way for any value of the coupling constant $\alpha$.
In~\cite{AGG}, the
moment map~$\mu_\alpha$ is computed, and it is shown that its zero locus corresponds to
solutions of~\eqref{eq:CKYM0}. The coupling between the metric and
the connection occurs as a direct consequence of the fact that
the extension~\eqref{eq:Ext-Lie-groupsintro} defining the extended gauge group is non-trivial.

It is worth pointing out that extended gauge groups feature in
the paper by Bourguignon--Lawson \cite{bourguignon-lawson}, where they are referred to as enlarged gauge groups. In particular they consider the {\em enlarged gauge group} of a principal bundle $P$ over a compact Riemannian manifold $M$ given~by
\[1 \to \GGG_P \to \widetilde{\GGG}_P \to \III_M \to 1,\]
where $\III_M$ is the group of isometries of $M$ when the dimension of $M$ is different from $4$ and the conformal group when the dimension is $4$. As mentioned in~\cite{bourguignon-lawson}, a connection on $P$ determines a~splitting of the sequence of vector spaces obtained by differentiating the above extension at the identity. This fact is also true in our set-up and plays a crucial role in the computation in~\cite{AGG} of the moment map for the action of $\widetilde{\GGG}$ on $\PPP$.

It turns out that equations \eqref{eq:CKYM0}
decouple on a compact Riemann surface, due to the vanishing of the first Pontryagin term $\tr F_H^2$, and so in this case the solution to the problem reduces to a combination of the uniformization theorem for Riemann surfaces and the theorem of Narasimhan and Seshadri~\cite{D2,NS}. For an arbitrary higher-dimensional manifold $M$, determining whether~\eqref{eq:CKYM0} admits solutions is a difficult problem, since in this case these equations are a~system of coupled fourth-order fully non-linear partial differential equations.

Despite this, a large class of examples was found in~\cite{AGG} for small $\alpha$, by perturbing constant scalar curvature K\"ahler metrics and Hermite--Yang--Mills connections. More concrete and interesting solutions over a polarised threefold --- that does not admit any constant scalar curvature K\"ahler metric
\textemdash{} were obtained by Keller and T{\o}nnesen--Friedman \cite{KellerTonnesen}. Garcia-Fernandez and Tipler~\cite{GFT} added new examples to this short list, by simultaneous deformation of the complex structures of $M$ and $E$.

However the problem of finding a general existence theorem for the K\"ahler--Yang--Mills equations remains pretty much open. In \cite{AGG}, one of
the main motivations to study these equations was to find an analytic approach to the
algebraic geometric problem of constructing a moduli space classifying pairs $(M,E)$ consisting of a complex projective variety and a holomorphic vector bundle. The stability condition needed for this should naturally be the one solving the K\"ahler--Yang--Mills
equations. In \cite{AGG}, obstructions for the existence of
solutions to \eqref{eq:CKYM0} were studied, generalizing the {\em Futaki invariant}, the {\em Mabuchi
$K$-energy} and {\em geodesic stability} that
appear in the constant scalar curvature theory. The natural conjecture proposed in~\cite{AGG} is that the existence of solutions
of the K\"ahler--Yang--Mills equations is equivalent to geodesic stability.

To test the above conjecture and provide a new class of interesting examples, inspired by~\cite{G1}, a series of papers \cite{AGG2,AGG3,AGGP,AGGPY2,AGGPY} have considered the study of dimensional reduction techniques. The simplest situation considered is given by
the dimensional reduction of the K\"ahler--Yang--Mills equations
from $M = X \times \PP^1$
to a compact Riemann surface $X$ of genus $g(X)$. Here $\PP^1$ is the Riemann sphere. In this case, we consider $\SU(2)$ acting on $M$, trivially on $X$ and in the standard way on $\PP^1$. We take a holomorphic line bundle $L$ over $X$, and a holomorphic global section $\phi$ of $L$. The pair $(L,\phi)$ defines an $\SU(2)$-equivariant holomorphic rank 2 vector bundle over $ X \times \PP^1$ in a canonical way. One can show \cite{AGG2} that an $\SU(2)$-invariant solution to the K\"ahler--Yang--Mills equations on the bundle $E$ over $X\times \PP^1$ is equivalent to having a solution of the equations
\begin{gather}\label{eq:gravvortexeq}
i\Lambda_g F_h + \frac{1}{2}\big(|\phi|_h^2-\tau\big) = 0,\qquad
S_g + \alpha(\Delta_g + \tau) \big(|\phi|_h^2 -\tau\big) = c
\end{gather}
for a K\"ahler metric $g$ on $X$ and a Hermitian metric $h$ on $L$. Here, $F_h$ is the curvature of the Chern connection on $L$ defined by $h$, $|\phi|_h$ is the pointwise norm of $\phi$ with respect to $h$, $S_g$ is the scalar curvature of $g$, and $\Delta_g$ is the Laplacian of the metric on the surface acting on functions. The constant $c\in\RR$ is topological, and it can be obtained by integrating~\eqref{eq:gravvortexeq} over $X$, and $\tau$ is a real parameter.

Equations \eqref{eq:gravvortexeq} are referred as the {\em gravitating vortex equations} since, in fact, the first equation in \eqref{eq:gravvortexeq} is the well-known vortex equation of the abelian Higgs model, whose solutions are called {\em vortices}, and have been extensively studied in the literature in the case of compact Riemann surfaces \cite{Brad,G1,G2,Noguchi} after the seminal work of A.~Jaffe and C.~Taubes~\cite{Jaffe-Taubes,Taubes1} on the Euclidean plane.
In particular, the proof given by S.~Bradlow \cite{Brad} is based on the fact that the vortex equation can be reduced to the Kazdan--Warner equation \cite{KaWa} --- an equation that plays a~prominent role in the problem studied by Bourguignon--Ezin \cite{bourguignon-ezin}.

It turns out that, when $c = 0$ and $d=c_1(L)>0$, $X$ is constrained to be the Riemann sphere and the gravitating vortex equations have a physical interpretation, as they are equivalent to the Einstein--Bogomol'nyi equations on a Riemann surface \cite{Yang1992,Yang1994CMP}. Solutions of the Einstein--Bogomol'nyi equations are known in the physics literature as
{\em Nielsen--Olesen cosmic strings}~\cite{NielsenOlesen}, and describe a special class of solutions of the abelian Higgs model coupled with gravity in four dimensions \cite{ComtetGibbons,Linet,Linet2}. Unlike the cases of genus $g(X)\geq 1$, in
genus $g(X)=0$, new phenomena, not appearing in the classical
situation of constant curvature metrics on a surface, arise. Namely, there exist obstructions to the existence of solutions of \eqref{eq:gravvortexeq},
illustrating the fact that, for $g(X)=0$, the problem of existence of solutions is comparatively closer to the more
sophisticated problem of Calabi on the existence of K\"ahler--Einstein
metrics, where algebro-geometric stability obstructions appear on
compact K\"ahler manifolds with $c_1>0$.

After presenting the K\"ahler--Yang--Mills equations in Section \ref{kym}, and reviewing in the theorems on the existence of solution to the gravitating vortex equations in Section~\ref{gve}, in Section \ref{non-Abelian-vortices}, we ponder on the existence of solutions for a non-abelian version of the gravitating vortex equations obtained also by dimensional reduction methods from the K\"ahler--Yang--Mills equations.

\section{The K\"ahler--Yang--Mills equations}\label{kym}

In this section, we briefly explain some basic facts
from~\cite{AGG} about the K\"ahler--Yang--Mills equations,
with emphasis on their symplectic interpretation.
Throughout this section, manifolds, bundles, metrics, and
similar objects are of class $C^\infty$. Let $M$ be a compact
symplectic manifold of dimension $2n$, with symplectic form $\omega$
and volume form $\vol_\omega =\frac{1}{n!}\omega^n$. Fix a complex
vector bundle $\pi\colon E\to M$ of rank $r$, and a Hermitian metric~$H$ on~$E$. Consider the positive definite inner product
\[
-\tr\colon\ \mathfrak{u}(r)\times\mathfrak{u}(r)\lto\RR
\]
on $\mathfrak{u}(r)$. Being invariant under the adjoint U$(r)$-action,
it induces a metric on the (adjoint) bundle $\ad E_H$ of
skew-Hermitian endomorphisms of $(E,H)$. Let $\Omega^k$ and
$\Omega^{k}(V)$ denote the spaces of (smooth) $k$-forms and $V$-valued
$k$-forms on $M$, respectively, for any vector bundle $V$ over
$M$. Then, the metric on $\ad E_H$ extends to a pairing on the space
$\Omega^\bullet(\ad E_H)$,
\[
\Omega^p(\ad E_H)\times\Omega^q(\ad E_H)\lto\Omega^{p+q},
\]
that will be denoted simply $-\tr a_p\wedge a_q$, for
$a_j\in\Omega^j(\ad E_H)$, $j=p,q$. An almost complex structure on $M$
compatible with $\omega$ determines a metric on $M$ and an operator
\begin{equation}
\label{eq:Lambda}
\Lambda\colon\ \Omega^{p,q}\lto\Omega^{p-1,q-1}
\end{equation}
acting on the space $\Omega^{p,q}$ of smooth $(p,q)$-forms, given by the adjoint of the Lefschetz operator~{$\Omega^{p-1,q-1}\to\Omega^{p,q}\colon\gamma\mapsto\gamma\wedge \omega$}. It can be seen that $\Lambda$ is symplectic, that is, it does not depend on the choice of almost complex structure on $M$. Its
linear extension to adjoint-bundle valued forms will also be denoted
$\Lambda\colon\Omega^{p,q}(\ad E_H)\to\Omega^{p-1,q-1}(\ad E_H)$.

Let $\JJJ$ and $\AAA$ be the spaces of almost complex structures on $M$
compatible with $\omega$ and unitary connections on $(E,H)$,
respectively; their respective elements will usually be denoted~$J$
and $A$. We will explain now how the K\"ahler--Yang--Mills equations
arise naturally in the construction of the symplectic quotient of a
subspace $\PPP\subset\JJJ\times\AAA$ of `integrable pairs'.

The group of symmetries of this theory is the \emph{extended gauge group} $\widetilde{\GGG}$. Let $E_H$ be the principal $\U(r)$-bundle of unitary frames of $(E,H)$. Then, $\widetilde{\GGG}$ is the group of automorphisms of~$E_H$ which cover elements of the group $\HHH$ of Hamiltonian symplectomorphisms of $(M,\omega)$. There is a canonical short exact sequence of Lie groups
\begin{equation}
\label{eq:coupling-term-moment-map-1}
 1\to \GGG \lra{} \widetilde{\GGG} \lra{\pr} \HHH \to 1,
\end{equation}
where $\pr$ maps each $g\in\widetilde{\GGG}$ into the Hamiltonian symplectomorphism
$\pr(g)\in\HHH$ that it covers, and so its kernel $\GGG$ is the unitary
gauge group of $(E,H)$, that is, the normal subgroup of $\widetilde{\GGG}$
consisting of unitary automorphisms covering the identity map on $M$.

There are $\widetilde{\GGG}$-actions on $\JJJ$ and $\AAA$, which, combined, give an
action on the product $\JJJ\times\AAA$,
\[
 g(J,A)=(\pr(g)J, gA).
\]
Here, $\pr(g)J$ denotes the push-forward of $J$ by $\pr(g)$. To
define the $\widetilde{\GGG}$-action on $\AAA$, we view the elements of $\AAA$ as
$G$-equivariant splittings $A\colon TE_H\to VE_H$ of the short exact
sequence
\[
 0 \to VE_H \lto TE_H\lto \pi^*TM \to 0,
\]
where $VE_H\subset TE_H$ is the vertical bundle on $E_H$.
Then, the $\widetilde{\GGG}$-action on $\AAA$ is $g A \defeq g\circ A \circ g^{-1}$, where $g\colon TE\to TE$ denotes the infinitesimal action on the right-hand side.

For each unitary connection $A$, we write $A^\perp y$ for the
corresponding horizontal lift of a vector field $y$ on $M$ to a vector
field on $E_H$. Then, each $A\in\AAA$ determines a vector-space
splitting of the Lie-algebra short exact sequence
\[
0\to\LieG\lto\LieX\lra{\pr}\LieH\to 0
\]
associated to~\eqref{eq:coupling-term-moment-map-1}, because $A^\perp\eta\in\LieX$ for all $\eta\in\LieH$. Note also that the equation
\[\eta_\varphi\lrcorner\omega={\rm d}\varphi\]
determines an isomorphism between the space $\LieH$ of Hamiltonian
vector fields on~$M$ and the space $C_0^\infty(M,\omega)$ of smooth
functions $\varphi$ such that $\int_M\varphi\vol_\omega=0$, where
$\vol_\omega\defeq\frac{1}{n!}\omega^n$.

The spaces $\JJJ$ and $\AAA$ have \smash{$\widetilde{\GGG}$}-invariant symplectic structures
$\omega_\JJJ$ and $\omega_\AAA$ induced by $\omega$, that, combined,
define a symplectic form on $\JJJ\times\AAA$, for each non-zero real
constant $\alpha$, given~by
\begin{gather}
\label{eq:Sympfamily}
\omega_\alpha = \omega_\JJJ + \frac{4 \alpha}{(n-1)!} \omega_\AAA.
\end{gather}
The following result provides the starting point for the theory of the
K\"ahler--Yang--Mills equations. This result builds on the moment map
interpretation of the constant scalar curvature equation for a
K\"ahler metric, due to Fujiki \cite{Fujiki} and Donaldson \cite{D1}, and the classical result of Atiyah and Bott \cite{AB}.

\begin{Proposition}[{\cite{AGG}}]
\label{prop:momentmap-pairs}
The $\widetilde{\GGG}$-action on $(\JJJ\times \AAA,\omega_\alpha)$ is Hamiltonian, with $\widetilde{\GGG}$-equivariant moment map \smash{$\mu_{\alpha}\colon \JJJ\times \AAA\to\bigl(\LieX\bigr)^*$} given by
\begin{align*}
\langle \mu_{\alpha}(J,A),\zeta\rangle ={}& 4{\rm i}\alpha\int_M \tr A\zeta\wedge({\rm i}\Lambda F_A-\lambda \Id)\vol_{\omega}\\
&- \int_M \varphi\(S_J - \alpha \Lambda^2 \tr F_A \wedge F_A - 4{\rm i}\lambda \alpha \Lambda \tr F_A\)\vol_{\omega}
\end{align*}
for any $\zeta\in\LieX$ covering $\eta_\varphi \in \LieH$, with $\varphi\in C_0^\infty(M,\omega)$.
\end{Proposition}

Here, $F_A$ is the curvature of $A$, $\lambda\in\RR$ is determined by
the topology of the bundle and the cohomology class $[\omega]\in H^2(M,\RR)$, and $S_J$ is the Hermitian scalar curvature of $J$. Explicitly,
\[
F_A = - A\bigl[A^\perp\cdot,A^\perp\cdot\bigr]\in\Omega^2(\ad E_H),
\qquad
\lambda=\frac{2\pi nc_1(E)\cdot[\omega]^{n-1}}{r[\omega]^n},
\]
with the convention $2\pi c_1(E)=[{\rm i}\tr F_A]$. A key observation in
\cite{AGG,GF} is that the space $\JJJ \times \AAA$ has a (formally
integrable) complex structure $\mathbf{I}$ preserved by the
\smash{$\widetilde{\GGG}$}-action, given by
\[
\mathbf{I}_{\mid(J,A)}(\gamma,a) = (J\gamma,-a(J \cdot))\qquad
\text{ for } (\gamma,a) \in T_J\JJJ \oplus T_A\AAA.
\]
For positive $\alpha$, $\mathbf{I}$ is compatible with the family of
symplectic structures \eqref{eq:Sympfamily}, and so it defines
K\"ahler structures on $\JJJ\times\AAA$. The condition $\alpha>0$ will
be assumed in the sequel.

Suppose now that there exist K\"{a}hler structures on $M$ with
K\"{a}hler form $\omega$. This means the subspace $\JJJ^i \subset \JJJ$
of integrable almost complex structures compatible with $\omega$ is
not empty. For each $J\in \JJJ^i$, let $\AAA^{1,1}_J\subset\AAA$ be the
subspace of connections $A$ with $F_A \in \Omega_J^{1,1}(\ad E_H)$,
where~$\Omega^{p,q}_J$ is the space of $(p,q)$-forms with respect to
$J$. Then the space of \emph{integrable pairs}
\[
 \PPP\subset \JJJ\times \AAA,
\]
consisting of elements $(J,A)$ with \smash{$J\in\JJJ^i$} and $A\in\AAA^{1,1}_J$,
is a $\widetilde{\GGG}$-invariant (possibly singular) K\"ahler submanifold. The
zero locus of the induced moment map $\mu_\alpha$ for the \smash{$\widetilde{\GGG}$}-action
on $\PPP$ corresponds precisely to the solutions of the (coupled)
\emph{K\"ahler--Yang--Mills equations}
\begin{gather}\label{eq:CKYM1}
{\rm i} \Lambda F_A = \lambda \Id,\qquad
S_J - \alpha \Lambda^2 \tr F_A \wedge F_A = c.
\end{gather}
Here, $S_J$ is the scalar curvature of the metric
$g_J=\omega(\cdot,J\cdot)$ and the constant $c\in\RR$ depends on~$\alpha$, the cohomology class of $\omega$ and the topology of $M$ and
$E$ (see \cite[Section 2]{AGG}).

One can express the K\"ahler--Yang--Mills equations from an alternative point of view in which we fix a compact complex manifold $X$ of dimension $n$, a K\"ahler class $\Omega \in
H^{1,1}(X)$ and a~holomorphic vector bundle $E$ over $X$. Then these
equations, for a fixed constant parameter~$\alpha\in\RR$, are
\begin{gather}\label{eq:CKYM2}
{\rm i}\Lambda_\omega F_H =\lambda\Id,\qquad
S_\omega-\alpha\Lambda_\omega^2\tr F_H\wedge F_H =c,
\end{gather}
where the unknowns are a K\"ahler metric on $X$ with K\"ahler form
$\omega$ in $\Omega$, and a Hermitian metric~$H$ on $E$. In this case,
$F_H$ is the curvature of the Chern connection $A_H$ of $H$ on $E$,
and~$S_\omega$ is the scalar curvature of the K\"ahler metric. Note
that the operator in~\eqref{eq:Lambda} depends on~$\omega$, and the
constant $c\in\RR$ depends on $\alpha$, $\Omega$ and the topology of
$X$ and $E$.

\section{The gravitating vortex equations}\label{gve}

Let $X$ be a compact connected Riemann surface of arbitrary genus. Let $L$ be a holomorphic line bundle over $X$ and $\phi \in H^0(X,L)$ a holomorphic section of $L$. We fix a parameter $0 < \tau \in \RR$, a~coupling constant $\alpha \in \RR$, and a real parameter $c$.

The \emph{gravitating vortex equations}, for a K\"ahler metric on $X$ with K\"ahler form $\omega$ and a~Hermitian metric $h$ on $L$, are
\begin{gather}\label{eq:gravvortexeq1}
{\rm i}\Lambda_\omega F_h + \frac{1}{2}\big(|\phi|_h^2-\tau\big) = 0,\qquad
S_\omega + \alpha(\Delta_\omega + \tau) \big(|\phi|_h^2 -\tau\big) = c.
\end{gather}
Here, $S_\omega$ is the scalar curvature of $\omega$,~$F_h$ stands for the curvature of the Chern connection of~$h$,~$|\phi|_h^2$ is the smooth function on $X$ given by the norm-square of $\phi$ with respect to $h$ and $\Delta_\omega$ is the Laplace operator for the metric $\omega$, defined by
\[
\Delta_Xf = 2{\rm i} \Lambda_\omega \dbar \partial f \qquad \textrm{ for } f \in C^\infty(\Sigma).
\]

We will now show how to derive the gravitating vortex equations
\eqref{eq:gravvortexeq1} as a dimensional reduction of the K\"ahler--Yang--Mills
equations \eqref{eq:CKYM2}. To do this, we associate to $(X,L,\phi)$ a~rank~2 holomorphic vector bundle $E$ over $X\times \PP^1$.
This is given as an extension
\begin{equation}\label{eq:bundleE}
0 \to p^*L \lto E \lto q^*\mathcal{O}_{\PP^1}(2) \to 0,
\end{equation}
where $p$ and $q$ are the projections from $X\times \PP^1$ to $X$ and $\PP^1$ respectively. By $\mathcal{O}_{\PP^1}(2)$ we denote as usual the holomorphic line bundle with Chern class $2$ on $\PP^1$, isomorphic to the holomorphic tangent bundle of $\PP^1$. Extensions as above are parametrised by
\[
H^1(X,p^*L \otimes q^*\cO_{\PP^1}(-2)) \cong H^0(X,L) \otimes H^1\bigl(\PP^1,\mathcal{O}_{\PP^1}(-2)\bigr) \cong H^0(X,L),
\]
and we choose $E$ to be the extension determined by $\phi$.
 Let $\SU(2)$ act trivially on $X$, and in the standard way on
$\PP^1 \cong \SU(2)/\U(1)$. This action can be lifted to trivial actions on~$E$ and $p^*L$ and the standard action on $\mathcal{O}_{\PP^1}(2)$. Since the induced actions on $H^0(X,L)$ and \smash{$H^1\big(\PP^1,\mathcal{O}_{\PP^1}(-2)\big) \cong H^0\big(\PP^1,\mathcal{O}_{\PP^1}\big)^* \cong \CC$} are trivial, $E$ is an $\SU(2)$-equivariant holomorphic vector bundle over $X\times \PP^1$.

For $\tau \in \RR_{>0}$, consider the $\SU(2)$-invariant K\"ahler metric on $X\times \PP^1$ whose K\"ahler form is~$\omega_\tau = p^*\omega + \frac{4}{\tau} q^*\omega_{\rm FS}$, where $\omega$ is a K\"ahler form on $X$ and $\omega_{\rm FS}$ is the Fubini--Study metric on~$\PP^1$, given in homogeneous coordinates by
\[
\omega_{\rm FS} = \frac{{\rm i} {\rm d}z \wedge {\rm d}\overline z}{\big(1+|z|^2 \big)^2}
\]
and such that $\int_{\PP^1}\omega_{\rm FS} = 2\pi$. Assuming that the coupling constants $\alpha$ in
\eqref{eq:gravvortexeq1}
and \eqref{eq:CKYM2} coincide, we have the following \cite{AGG2}.

\begin{Proposition}\label{prop:dimred}
The triple $(X,L,\phi)$ admits a solution $(\omega,h)$ of the gravitating vortex equations~\eqref{eq:gravvortexeq1} with parameter $\tau$ if and only if $(X\times \PP^1,E)$ admits an $\SU(2)$-invariant solution of the K\"ahler--Yang--Mills equations \eqref{eq:CKYM2} with K\"ahler form $\omega_\tau= p^*\omega + \frac{4}{\tau} q^*\omega_{\rm FS}$.
\end{Proposition}

For a fixed K\"ahler metric $\omega$, the first equation in
\eqref{eq:gravvortexeq1} corresponds to the abelian vortex equation
\begin{equation}\label{eq:vortexeq}
{\rm i}\Lambda_\omega F_h + \frac{1}{2}\big(|\phi|_h^2-\tau\big) = 0,
\end{equation}
for a Hermitian metric $h$ on $L$. In \cite{Brad,G1,G2,Noguchi}, Noguchi, Bradlow
and the author gave independently and with different methods a complete characterisation of the existence of \emph{abelian vortices} on a compact Riemann surface, that is, of solutions of equations \eqref{eq:vortexeq}.

\begin{Theorem}[\cite{Brad,G1,G2}]\label{th:B-GP}
Assume that $\phi$ is not identically zero. For every fixed K\"ahler form~$\omega$, there exists a unique solution $h$ of the vortex equations \eqref{eq:vortexeq} if and only if
\begin{equation}\label{eq:ineq}
c_1(L) < \frac{\tau \Vol_\omega(X)}{4\pi}.
\end{equation}
\end{Theorem}

Inspired by work of Witten \cite{Witten} and Taubes \cite{Taubes}, the method in \cite{G1} exploited the dimensional reduction of the Hermitian--Yang--Mills equations from four to two dimensions, combined with the theorem of Donaldson, Uhlenbeck and Yau \cite{D3,UY,UY+}.

The constant $c \in \RR$ is topological, and is explicitly given by
\begin{equation}\label{eq:constantc}
c = \frac{2\pi(\chi(X) - 2\alpha\tau c_1(L))}{\Vol_\omega(X)},
\end{equation}
as can be deduced by integrating the equations.
The gravitating vortex equations for $\phi = 0$, are equivalent to the condition that
$\omega$ be a constant scalar curvature K\"ahler metric on $X$ and~$h$ be
a Hermite--Einstein metric on $L$. By the uniformisation theorem for Riemann
surfaces, the existence of these `trivial
solutions' reduces by Hodge theory to the condition ${c_1(L) = \tau\Vol_\omega(X)/4\pi}$.

Excluding this trivial case, the sign of $c$ plays an important role in the
existence problem for the gravitating vortex equations.
The dependence of the gravitating vortex equations~\eqref{eq:gravvortexeq1} on the topological constant $c$ is better observed using a K\"ahler--Einstein type formulation. Using that $X$ is compact, \eqref{eq:gravvortexeq1} reduces to a second-order system of PDE. To see this, we fix a constant scalar curvature metric $\omega_0$ on $X$ and the unique Hermitian metric $h_0$ on $L$ with constant $\Lambda_{\omega_0} F_{h_0}$, and apply a conformal change to $h$ while changing $\omega$ within its K\"ahler class. Equations \eqref{eq:gravvortexeq1} for~$\omega = \omega_0 + dd^c v$, $h={\rm e}^{2f}h_0$, with $v,f\in C^\infty(\Sigma)$, are equivalent to
the following semi-linear system of partial differential equations (cf.\ \cite[Lemma~4.3]{AGG2})
\begin{gather}
\Delta f + \frac{1}{2}\big({\rm e}^{2f}|\phi|^2-\tau\big){\rm e}^{4\alpha \tau f - 2 \alpha {\rm e}^{2f}|\phi|^2 - 2 c v} = - c_1(L),\nonumber\\
\Delta v + {\rm e}^{4\alpha \tau f - 2 \alpha {\rm e}^{2f}|\phi|^2 - 2 c v} = 1.\label{eq:KWtype0}
\end{gather}
Here, $\Delta$ is the Laplacian of the fixed metric $\omega_0$, normalised to have volume $2\pi$ and $|\phi|$ is the pointwise norm with respect to the fixed metric $h_0$ on $L$. Note that $\omega = (1- \Delta v) \omega_0$ implies~${1 - \Delta v >0}$, which is compatible with the last equation in \eqref{eq:KWtype0}.

For $c\geq 0$, the existence of gravitating vortices forces the topology of the surface to be that of the $2$-sphere, because $c_1(L)>0$ implies $\chi(\Sigma) > 0$ by~\eqref{eq:constantc}.
When $c$ in~\eqref{eq:constantc} is zero, the gravitating vortex
equations~\eqref{eq:gravvortexeq1} turn out to be a system of partial
differential equations that have been extensively studied in the
physics literature, known as the \emph{Einstein--Bogomol'nyi equations}.
As observed by Yang~\cite[Section 1.2.1]{Yang}, the existence of solutions in this situation
with $\alpha>0$
constrains the topology of $X$ to be the complex projective
line (or 2-sphere) $\PP^1$, since $c=0$ if and only if
\[
\chi(X) = 2\alpha\tau c_1(L).
\]
We are assuming that $\tau > 0$ and $c_1(L) > 0$.

In the case $c=0$, for $L=\mathcal{O}_{\PP^1}(N)$ and ${\rm e}^{2u} = 1- \Delta v$ the system \eqref{eq:KWtype0} reduces to a single partial differential equation
\begin{gather}\label{eq:single}
\Delta f + \frac{1}{2}{\rm e}^{2u}\big({\rm e}^{2f}|\phi|^2-\tau\big) = - N,
\end{gather}
for a function $f\in C^\infty\big(\PP^1\big)$, where
\[
u=2\alpha\tau f-\alpha {\rm e}^{2f}|\phi|^2+c',
\]
and $c'$ is a real constant that can be chosen at will. By studying the {\em Liouville type equation}~\eqref{eq:single} on $\PP^1$, Yang \cite{Yang,Yang3} proved the
existence of solutions of the Einstein--Bogomol'nyi equations under
certain numerical conditions on the zeros of $\phi$, to which he refers as a ``technical restriction''~\cite[Section 1.3]{Yang}. It turns out
that these conditions have a precise algebro-geometric meaning in the
context of Mumford's geometric invariant theory (GIT)~\cite{MFK}, as a
consequence of the following result.

\begin{Proposition}[{\cite[Chapter~4, Proposition~4.1]{MFK}}]\label{prop:GIT}
Consider the space of effective divisors on~$\PP^1$ with its canonical linearised $\SL(2,\CC)$-action. Let
\smash{$D=\sum_j n_jp_j$} be
an effective divisor, with finitely many different points $p_j\in\PP^1$
and integers $n_j>0$ such that~$N~=~\sum_j n_j$. Then
\begin{enumerate}\itemsep=0pt
\item[\textup{(1)}] $D$ is stable if and only if $n_j < \frac{1}{2}N$ for all $j$;
\item[\textup{(2)}] $D$ is strictly polystable if and only if $D=\frac{1}{2}Np_1 + \frac{1}{2}Np_2$, where $p_1 \neq p_2$ and $N$ is even;
\item[\textup{(3)}] $D$ is unstable if and only if there exists $p_j
 \in D$ such that $n_j>\frac{1}{2}N$.
\end{enumerate}
\end{Proposition}

Using Proposition~\ref{prop:GIT}, Yang's existence theorem has the
following reformulation, where ``GIT polystable'' means either
conditions (1) or (2) of Proposition~\ref{prop:GIT} are satisfied, and
\[
D = \sum_j n_j p_j
\]
is the effective divisor on $\PP^1$ corresponding to a pair
$(L,\phi)$, with $N = \sum_j n_j = c_1(L)$.

\begin{Theorem}[Yang's existence theorem]\label{th:Yang}
Assume that $\alpha > 0$ and that \eqref{eq:ineq} holds. Then, there exists a solution of the Einstein--Bogomol'nyi equations on $\big(\PP^1,L,\phi\big)$ if $D$ is GIT polystable for the
linearised $\operatorname{SL}(2,\CC)$-action on the space of effective divisors.
\end{Theorem}

The converse to Theorem~\ref{th:Yang} is given in \cite{AGGP,AGGPY2}.

\begin{Theorem}\label{th:Yangconjectureintro}
If $\big(\PP^1,L,\phi\big)$ admits a solution of the gravitating vortex
equations with $\alpha > 0$, then \eqref{eq:ineq} holds and the
divisor $D$ is polystable for the $\textup{SL}(2,\CC)$-action.
\end{Theorem}

\begin{Remark}
 Notice that this theorem is more general than being a converse to Theorem~\ref{th:Yang} since it does not assume
 that $c=0$, and deals with the general gravitating vortex equations~\eqref{eq:gravvortexeq1} and not just with the Einstein--Bogomol'nyi equation.
\end{Remark}

Combining now Theorems~\ref{th:Yang}
and~\ref{th:Yangconjectureintro}, we obtain a
correspondence theorem for the Einstein--Bogomol'nyi equations.

\begin{Theorem}\label{th:HK}
A triple $\big(\PP^1,L,\phi\big)$ with $\phi\neq 0$
admits a solution of the Einstein--Bogomol'nyi equations with $\alpha > 0$ if and only
\eqref{eq:ineq} holds and the divisor $D$ is polystable for the
$\textup{SL}(2,\CC)$-action.
\end{Theorem}

Another result, conjectured by Yang and proved in \cite{AGGP} is the
following.

\begin{Theorem}
There is no solution of the Einstein--Bogomol'nyi equations for $N$
strings superimposed at a single point, that is, when $D=Np$.
\end{Theorem}

An existence theorem for the gravitating vortex equations~\eqref{eq:gravvortexeq1} with $c>0$ for a triple $\big(\PP^1,L,\phi\big)$ with $\phi\neq 0$, similar to Theorem \ref{th:HK}, is obtained combining Theorem \ref{th:Yangconjectureintro} with the converse direction in this situation, proved by Garcia-Fernandez--Pingali--Yao \cite{GPY}.

In genus $g(X)=1$, the gravitating vortex equations~\eqref{eq:gravvortexeq1}
(with $\phi\neq 0$) always have a~solution in the weak coupling limit
$0< \alpha \ll 1$ (see~\cite[Theorem~4.1]{AGG2} for a
precise formulation), and it is an interesting open problem to find
effective bounds for $\alpha$ for which~\eqref{eq:gravvortexeq1} admit
solutions.

Paper \cite{AGGP} deals also with the existence theorem
of solutions of~\eqref{eq:gravvortexeq1} for surfaces of genus~$g\geq 2$, for which one has the following.

\begin{Theorem}
\label{thm:higher-genus.intro}
Let $X$ be a compact Riemann surface of genus $g\geq 2$, and $L$
a holomorphic line bundle over $X$ of degree $N>0$ equipped with a
holomorphic section $\phi \neq 0$. Let $\tau$ be a real constant such
that $0<N<\tau/2$. Define
\begin{equation}
\label{eq:critical-alpha.intro}
\alpha_*\defeq\frac{2g-2}{2\tau(\tau/2-N)}>0.
\end{equation}
Then, the set of $\alpha$ for which~\eqref{eq:gravvortexeq} admits
smooth solutions of volume $2\pi$ is open and contains the closed
interval $[0,\alpha_*]$. Furthermore, the solution is unique for
$\alpha\in[0,\alpha_*]$.
\end{Theorem}

This shall be compared with the classical uniformization theorem which establishes
that a~compact Riemann surface admits a metric of constant curvature with fixed
volume, unique up to biholomorphisms. The proof of Theorem~\ref{thm:higher-genus.intro} involves the
continuity method, where openness is proven using the moment-map
interpretation, while closedness needs
\emph{a priori} estimates as usual.
The hardest part is the $C^0$ estimate, and in fact it is
for this estimate that the value of $\alpha$ should not be too
large. With these estimates at hand, we prove uniqueness by adapting an argument by Bando and Mabuchi
in the K\"ahler--Einstein situation~\cite{BM}.
An interesting open question is to
see what the largest value of $\alpha$ is, for which solutions
exist. Notice that in the \emph{dissolving limit} $\tau \to N/2$ of the vortex we have $\phi \to 0$ (see \cite{G1}), and $\alpha^*$ in \eqref{eq:critical-alpha.intro} becomes arbitrarily large.

\section{Non-abelian gravitating vortices}\label{non-Abelian-vortices}

One can consider the dimensional reduction of the K\"ahler--Yang--Mills equations for higher rank $\SU(2)$-equivariant bundles on $X\times \PP^1$, where $X$ is a compact Riemann surface. In particular one can consider extensions of the form
\begin{equation}\label{eq:triple}
0 \to p^*E_1 \lto E \lto p^*E_2\otimes q^*\mathcal{O}_{\PP^1}(2) \to 0,
\end{equation}
where $E_1$ and $E_2$ are holomorphic vector bundles on $X$ and, as above, $p$ and $q$ are the projections from $X\times \PP^1$ to $X$ and $\PP^1$ respectively. Extensions of the form \eqref{eq:triple} are in one-to-one correspondence with triples
$T=(E_1,E_2,\phi)$, where $\phi$ is a sheaf homomorphism from $E_2$ to $E_1$, that is an element in $H^0(X,\Hom(E_2,E_1))$. As
\eqref{eq:bundleE}, these extensions define $\SU(2)$-equivariant (in fact, $\SL(2,\CC)$-equivariant) vector bundles over $X\times \PP^1$.

Given a triple $T=(E_1,E_2,\phi)$ over $X$ we can consider the
{\em gravitating coupled vortex equations} for a metric on $X$ with K\"ahler form $\omega$, and Hermitian metrics $h_1$ and $h_2$ on $E_1$ and $E_2$, respectively, given by
\begin{gather}
 {\rm i}\Lambda_\omega F_{h_1} + \phi\phi^* = \tau_1,\qquad
 {\rm i}\Lambda_\omega F_{h_2} - \phi^*\phi = \tau_2,\nonumber\\
 S_\omega + \alpha\Delta_\omega |\phi|^2
 -\alpha\big(\Lambda_\omega^2 \tr F_{h_1}^2+\Lambda_\omega^2 \tr F_{h_2}^2 +4\tau_1 \tr({\rm i}\Lambda_\omega F_{h_1})
 +4\tau_2 \tr({\rm i}\Lambda_\omega F_{h_2})\big) = c.\label{eq:non-Abelian-gve}
\end{gather}

See Section \ref{gve} for the notations. Here $\tau_1$ and $\tau_2$ are real parameters (of which a certain linear combination is linked to the Chern classes of $E_1$ and $E_2$), so that $\tau_1-\tau_2>0$, $\alpha\in \RR$ is a coupling constant, and $c\in \RR$ depends on $\alpha$, $\tau_1-\tau_2$, and the topology of $X$ and $E_1$ and $E_2$.
\begin{Remark}
When referring to the gravitating coupled vortex equations we are of course assuming that at least one of the two vector bundles has rank bigger than one. The case in which both vector bundles have rank one can be reduced to the study of the abelian gravitating vortex equation as shown in \cite{AGG}.
\end{Remark}

Let $\sigma:=\tau_1-\tau_2$, and consider now the $\SU(2)$-equivariant K\"ahler form
\[\omega_\sigma=\sigma p^*\omega + q^*\omega_{\rm FS},\] where $\omega$ is a K\"ahler form on $X$ and $\omega_{\rm FS}$ is the Fubini--Study metric on $\PP^1$ (see Section \ref{gve}).
Similarly to Proposition \ref{prop:dimred}, one has the following \cite{AGG3}.

\begin{Proposition}\label{prop:dimred-non-Abelian}
 Let $T=(E_1,E_2,\phi)$ be a triple over $X$. The pair $(X,T)$ admits a solution $(\omega,h_1,h_2)$ of the gravitating coupled
 vortex equations \eqref{eq:non-Abelian-gve} if and only if $(X \times \PP^1, E )$ admits an $\SU(2)$-invariant solution of the K\"ahler--Yang--Mills equations \eqref{eq:CKYM2} with K\"ahler form~$\omega_\sigma$.
\end{Proposition}

For a fixed K\"ahler metric $\omega$, the first two equations in
\eqref{eq:non-Abelian-gve} are the {\em coupled vortex equations} introduced in \cite{G3}, where it was shown that they
are dimensional reduction of the Hermite--Yang--Mills equations. An existence theorem for the coupled vortex equations was given
in \cite{BG} in terms of a certain notion of stability for the triple $T$ depending on the parameter $\sigma$. To define this concept, let $T = (E_{1},E_{2},\phi)$ and $T'~=~(E_{1}',E_{2}',\phi')$ be two triples on $X$.
A homomorphism from $T'$ to $T $ is a commutative diagram
\begin{displaymath}
 \begin{CD}
 E_2' @>\phi'>> E_1' \\
 @VVV @VVV \\
 E_2 @>\phi>> E_1,
 \end{CD}
\end{displaymath}
where the vertical arrows are holomorphic maps.
A triple $T'=(E_1',E_2',\phi')$ is a subtriple of~${T = (E_1,E_2,\phi)}$
if the sheaf homomorphims $E_1'\to E_1$ and $E_2'\to E_2$
are injective. A
subtriple~${T'\subset T}$ is called {\em proper} if
 $T'\neq 0 $ and $T'\neq T$.

For any $\sigma \in \RR$ the \emph{$\sigma$-degree} and
\emph{$\sigma$-slope} of $T$ are
defined to be{\samepage
\begin{align*}
 &\deg_{\sigma}(T)
 = \deg(E_{1}) + \deg(E_{2}) + \sigma
 \rk(E_{2}), \\
 &\mu_{\sigma}(T)
 =
 \frac{\deg_{\sigma}(T)}
 {\rk(E_{1})+\rk(E_{2})} = \mu(E_{1} \oplus E_{2}) +
 \sigma\frac{\rk(E_{2})}{\rk(E_{1})+
 \rk(E_{2})},
\end{align*}
where $\deg(E)$, $\rk(E)$ and $\mu(E)=\deg(E)/\rk(E)$ are the
degree, rank and slope of $E$, respectively.}

We say $T\! = (E_{1},E_{2},\phi)$ is
\emph{$\sigma$-stable} if $\mu_{\sigma}(T')
 < \mu_{\sigma}(T)$ for any proper subtriple $T' \!= (E_{1}',E_{2}',\phi')$.
We define \emph{$\sigma$-semi-stability} by replacing the above
strict inequality with a weak inequality. A~triple is called
\emph{$\sigma$-polystable} if it is the direct sum of $\sigma$-stable
triples of the same $\sigma$-slope.
We denote by
$ \mathcal{M}_\sigma
 = \mathcal{M}_\sigma(n_1,n_2,d_1,d_2)$
the moduli space of $\sigma$-polystable triples $T =
(E_{1},E_{2},\phi)$ which have $\rk(E_i)=n_i$ and $\deg(E_i) = d_i$ for
$i=1,2$.

There are certain necessary conditions in order for $\sigma$-semistable
triples to exist.
Let $\mu_i=d_i/n_i$ for $i=1,2$. We define
\begin{align*}
 &\sigma_m= \mu_1-\mu_2, \qquad
 \sigma_M = \left(1+ \frac{n_1+n_2}{|n_1 - n_2|}\right)(\mu_1 - \mu_2), \qquad
 n_1\neq n_2.
\end{align*}
\begin{Proposition}[{\cite[Theorem 6.1]{BG}}]
 \label{prop:alpha-range}
The moduli space $\mathcal{M}_\sigma(n_1,n_2,d_1,d_2)$ is a complex analytic
variety, which is projective when $\sigma$ is rational.
A necessary condition for $\mathcal{M}_\sigma(n_1,n_2,d_1,d_2)$
to be non-empty is
\begin{gather*}
0\leq \sigma_m \leq \sigma \leq \sigma_M \quad \text{if}\quad n_1\neq n_2,\qquad
0\leq \sigma_m \leq \sigma\quad \text{if}\quad n_1= n_2.
\end{gather*}
\end{Proposition}

The non-emptiness and topology of these moduli spaces have been studied in
\cite{BGG1,GH,GHS}. Triples play an important role in the study of the moduli space of Higgs bundles and character varieties
of the fundamental group of~$X$ \cite{BGG1,BGG2}.

The study of existence of solutions to the gravitating coupled vortex equations \eqref{eq:non-Abelian-gve} in the higher rank case is entirely open,
to the knowledge of the author. Solving this problem will most likely require new analytic and algebraic
tools and techniques, that may turn out to be very useful for the study of existence of solutions of the K\"ahler--Yang--Mills equations~\eqref{eq:CKYM1}. A particularly interesting situation, similarly to the abelian case should be the case in which~${X=\PP^1}$. In this situation, one may conjecture the following.

\begin{Conjecture} Let $T = (E_{1},E_{2},\phi)$ be a triple over $\PP^1$. The pair $(\PP^1,T)$ admits a solution to the gravitating coupled vortex equations \eqref{eq:non-Abelian-gve} if and only $T$ is $\sigma$-polystable and the point~${T\in \mathcal{M}_\sigma}$ is GIT polystable for the natural action of $\SL(2,\CC)$ on $\mathcal{M}_\sigma$ induced by the action of~$\SL(2,\CC)$ on~$\PP^1$.
\end{Conjecture}

\begin{Remark}
When $n_1=n_2=1$, $\sigma$-stability of the triple $T$ reduces to the condition~\eqref{eq:ineq}, where $L=E_1\otimes E_2^*$ and $\sigma$ is essentially the inverse of $\tau$. Then, the proof of this conjecture reduces to Theorem~\ref{th:HK}.
\end{Remark}

\subsection*{Acknowledgements}

The author thanks his co-authors on the various
subjects treated in this paper. These include: Luis \'Alvarez-C\'onsul, Steven Bradlow, Mario Garcia-Fernandez, Peter Gothen, Vamsi Pingali and Chengjian Yao. He also thanks Jean-Pierre Bourguignon for comments and corrections on the first draft of this paper, and
the IHES for its hospitality and support.
Partially supported by the Spanish Ministry of Science and
Innovation, through the ``Severo Ochoa Programme for Centres of Excellence in R\&D (CEX2019-000904-S)'' and PID2022-141387NB-C21.

\pdfbookmark[1]{References}{ref}
\LastPageEnding

\end{document}